\documentclass[10pt]{amsart}

\usepackage{tikz}
\usetikzlibrary{positioning}
\usetikzlibrary{arrows}
\usepackage{amssymb,amsmath,latexsym,graphicx}
\usepackage{charter,eucal}
\usepackage{amssymb}
\usepackage{amscd}
\usepackage[all,cmtip]{xy}
\usepackage{rotating,multicol,soul}
\usetikzlibrary{decorations.markings}
\usepackage{tcolorbox,upgreek,fontenc,mathtools}

\newcommand*\circled[1]{\tikz[baseline=(char.base)]{
            \node[shape=circle,draw,inner sep=2pt] (char) {#1};}}

\setul{}{1.1pt}

\newtheorem{theorem}{Theorem }[section]

\newtheorem{lemma}[theorem]{Lemma}
\newtheorem{observation}[theorem]{Observation}
\newtheorem{problem}[theorem]{Problem}
\newtheorem{remark}[theorem]{Remark}
\newtheorem{corollary}[theorem]{Corollary}
\newtheorem{proposition}[theorem]{Proposition}
\newtheorem{question}[theorem]{\textsc{Question}}

\newcommand{\bt}{\begin{theorem}}
\newcommand{\et}{\end{theorem}}
\newcommand{\bmt}{\begin{maintheorem}}
\newcommand{\emt}{\end{maintheorem}}
\newcommand{\bc}{\begin{corollary}}
\newcommand{\bl}{\begin{lemma}}
\newcommand{\ec}{\end{corollary}}
\newcommand{\el}{\end{lemma}}
\newcommand{\bo}{\begin{observation}}
\newcommand{\eo}{\end{observation}}
\newcommand{\bp}{\begin{proposition}}
\newcommand{\ep}{\end{proposition}}
\newcommand{\br}{\begin{remark}}
\newcommand{\er}{\end{remark}}
\newcommand{\bpr}{\begin{principle}}
\newcommand{\epr}{\end{principle}}
\newcommand{\bq}{\begin{question}}
\newcommand{\eq}{\end{question}}
\newcommand{\bpro}{\begin{problem}}
\newcommand{\epro}{\end{problem}}

\def\Aut{\mathrm{Aut}}
\def\I{\mathop{\mathrm{I}}}
\def\nI{\mathop{\not\bI}}

\def\I{\mathbf{I}}

\def\eop{\hspace*{\fill}$\blacksquare$}

\def\id{\mathrm{id}}

\newcommand{\F}{\mathbb{F}}
\newcommand{\bP}{\mathbb{P}}
\newcommand{\mC}{\mathcal{C}}

\newcommand{\mQ}{\mathcal{Q}}
\newcommand{\A}{\mathbb{A}}
\newcommand{\mA}{\mathcal{A}}
\newcommand{\mP}{\mathcal{P}}
\newcommand{\mL}{\mathcal{L}}
\newcommand{\mE}{\mathcal{E}}

\newcommand{\bI}{\mathbf{I}}

\newcommand{\mO}{\mathcal{O}}

\newcommand{\mI}{\mathcal{I}}

\newcommand{\mD}{\mathcal{D}}

\newcommand{\ob}{\mathrm{ob}}
\newcommand{\mB}{\mathcal{B}}

\newcommand{\mU}{\mathcal{U}}

\newcommand{\wt}{\widetilde}

\usetikzlibrary{matrix}
\usetikzlibrary{arrows}
\usepackage{pgf}
\usepackage{lscape}
\usepackage{listings}

\usepackage{enumitem}

\title{Synthetic projective lines, geometric closure and AB-sets}
\author{Koen Thas}
\thanks{}
\address{Ghent University, Department of Mathematics, Krijgslaan 281, S25, B-9000, Ghent, Belgium}
\email{koen.thas@gmail.com}\date{}

\begin{document}

\maketitle
\begin{abstract}
In this note, we introduce a new approach to abstract ``synthetic'' projective lines. We discuss various aspects of our approach, and compare these aspects 
with the classical one. A number of intriguing questions arise. Amongst these aspects, we discuss geometric closures, and introduce {\em automorphism-blocking sets}. We also construct a number of (counter) examples in infinite cases. 
\end{abstract}

\begin{tcolorbox}
\tableofcontents
\end{tcolorbox}

\section{Introduction}

There exists a long and rich history in defining the concept of a projective line, both over coordinatizing algebraic structures and combinatorially. 
The classical approach is to consider a field (or division ring) $k$, and define the projective line $\bP^1(k)$ over $k$ as the geometry of the vector lines in the $2$-dimensional (left or right, in case of a division ring) vector space $V(2,k)$. (One can see the line as the set $k \cup \{ \infty \}$, which comes with the natural $3$-transitive action of $\mathsf{PGL}_2(k)$, extended with field automorphisms.) This local viewpoint is a property we want to incorporate in our approach.

\subsection*{Some known approaches}

One first, naive synthetic approach is the rather ``flat'' incidence-geometrical way in which one sees a line as a set of points $\Omega$, which then comes with the automorphism 
group $\mathsf{Sym}(\Omega)$. This is obviously not what we want: their is no projective aspect whatsoever, and the automorphism group is boring. Buekenhout \cite{Buek}, under the influence of the school of Libois,  
took an entire different turn, and modeled projective lines after projective lines seen as subgeometries of Desarguesian projective spaces $\mathbb{P}^n(k)$ (with $k$ a field). More specifically, he molded the geometry by using specific automorphisms (which in the Desarguesian examples are induced by central collineations in the line stabilizers in $\mathsf{P\Upgamma L}_{n + 1}(k)$), as such creating an interesting geometric object which comes with an interesting automorphism group. We recall Buekenhout's definition. A {\em projective line} is a set $\mL$ of cardinality $\geq 3$,  satisfying the following axioms:
\begin{itemize}
\item[(a)]
Each $p, q \in \mL$, comes with a group $\lambda(p, q)$ whose elements are called {\em central collineations} with {\em center} $p$ and {\em axis} $q$. Note that $q$ need not be distinct from $p$. The group $\lambda(p, q)$ acts on $\mL$ as follows: it fixes $p$ and $q$, and if an element fixes a point $r \not\in \{ p, q \}$, then it is the identity. The action is faithful, and each $\lambda(p,q)$ is a subgroup of $\Aut(\mL)$.
\item[(b)]
If $p \ne q$, then $\lambda(p,q)$ and $\lambda(q,p)$ commute with each other.
\item[(c)]
For any $\gamma \in \lambda(r,s)$ and any $p, q$, we have $\gamma\lambda(p,q)\gamma^{-1} = \lambda(\gamma(p),\gamma(q))$.
\item[(d)]
The composite of two collineations with center $p$ (possibly with different axes) is again a collineation with center $p$, and dually the composite of two collineations with axis $q$ (possibly with different centers) is again a collineation with axis $q$. Thus we have two groups $\lambda(p) = \cup_q\lambda(p,q)$ and $\lambda(q) = \cup_p\lambda(p,q)$.
\end{itemize}

A projective line is {\em Desarguesian} if in addition $\lambda(p,q)$ acts (sharply) transitively on $\mL \setminus \{ p,q\}$. 
In this approach, any set is a projective line (by putting all $\lambda(p,q) = \{\id\}$). The natural examples are the lines $U$ in any (not necessarily Desarguesian) projective plane $\mP$: if $u$ and $v$ are points of $U$, then let $\lambda(u,v)$ be induced by the group of central collineations of $\mP$ with center $u$ and axis some line $V$ incident with $v$ where $V$ varies over these lines.  
Buekenhout's theory is just one of many quite similar approaches which endow a set with specific permutations modeled after lines in certain classes of projective planes. Still, it only models the automorphism groups of projective lines (in projective planes), and not the geometry. In our approach, projective lines will come with an embedding; this viewpoint includes a geometric theory, a relative angle and heavily restricted automorphism groups. 

One other way to see projective lines in the incidence-geometrical context is through Tits's theory of (split) BN-pairs, in this case of dimension $1$. Unfortunately this is essentially an approach to classical lines, and if $U$ is a line in a projective line $\mP$, then if $H$ is the group induced on $\Aut(\mP)_U$, the permutation group $(H,U)$ (where we see $U$ as a point set) is generically not reached by a (split) BN-pair of rank $1$, since the latter are always doubly transitive. 
We feel that this approach also lacks an underlying geometry, and in this paper we want to take quite a different view on the subject which entertains both more geometry and more general automorphism groups.

\begin{remark}{\rm 
One of the most interesting viewpoints is the algebro-geometric approach, which puts a Zariski topology on the set of prime ideals of a commutative ring (this 
defines the geometry), and then comes with a sheaf of rings. For our purposes, especially the geometrical (topological) aspect is very interesting, and the structure of the commutative ring endows the projective line with a highly nontrivial geometry. We give two examples. Consider the real projective  line $\mathsf{Proj}(\mathbb{R}[x, y])$. Its prime ideals (closed points) are generated by either polynomials of degree $1$ or irreducible polynomials of degree $2$. Those of degree $1$, such as $(x - a)$ gives a ``classical rational point'' ($(1,a)$) and those of the second type, such as $(x^2 + 1)$ ``become'' two rational points $(x + i)$ ($(1,i)$) and $(x - i)$ ($(1,-i)$) in the projective line $\mathsf{Proj}(\mathbb{C}[x, y])$, and $\mathsf{Gal}(\mathbb{C}/\mathbb{R})$ acts transitively on these ``conjugate points.'' {\em Second example}: consider the algebro-geometric projective line $L = \mathsf{Proj}({\F_q}[x,y])$, where $\F_q$ is any finite field. Its closed points are the ``usual'' points, but it contains much more information. Let $x$ be any of its closed points (which we regard as the point at infinity), and consider the affine line $L \setminus \{ x\} = \mathsf{Spec}(\F_q[x])$. Let $f(x) \in \F_q[x]$ be an irreducible polynomial of degree $m$. Then the point corresponding to the prime ideal $(f(x))$ defines $m$ different points in $\F_{q^m} \cong \F_q[x]/(f(x))$ on which the Galois group $\mathsf{Gal}(\F_{q^m}/\F_q)$ acts transitively. In some sense, the affine line $\mathsf{Spec}(\overline{\F_q}[x])$ contains the same information (and its Galois group captures the information of all the separate finite Galois groups).    

} 
\end{remark}

\subsection*{The present paper}

In this paper, a projective line comes with a full imbedding in an axiomatic projective plane, and morphisms between lines must take the embeddings into account. But there is more: if $U$ is a line in a projective plane $\mP$, then the actual projective line is determined by a point $x$ not incident with $U$, and that line consists of the data $\mL = (\mP, U, x)$. The points of the line $\mL$ are given by the lines incident with $x$, and automorphism must respect the structure of the plane, while fixing $U$ and $x$. We compare two approaches: 
the ``projective'' approach in which $U$ itself is defined as a line, and the aforementioned ``affine'' approach, which looks at the line from the viewpoint of $x$. In the Desarguesian case, those approaches are essentially the same, but in general they are very different. Generically, what one sees locally from the point $x$ depends on the location of $x$, and only when there is an automorphism of the plane which fixes $U$ while carrying $x$ to another point $x'$, the local geometries will be the same. On the other hand, it could be that nonisomorphic lines in this formalism define the exact same automorphism group as a permutation group, and that is different in the synthetic approached we mentioned (where the isomorphism class of the line is determined by the isomorphism class of the permutation group it defines). 

Our main viewpoint is essentially the following. If one wants to define abstract synthetic projective planes $\Upgamma$, one cannot use an embedding into a higher-dimensional projective space such as a $3$-space ($\circled{2}$), since the classical Veblen-Young result states that such as space is isomorphic to a space $\mathbb{P}^3(\dj)$ which is defined 
over a division ring $\dj$. But since planes  do not have this property, we {\em can} and will use embeddings of lines in axiomatic planes ($\circled{1}$) to do the job.  

\begin{equation}
\mbox{line}\ \xrightarrow[]{\circled{1}}\ \mbox{plane} \ \xrightarrow[]{\circled{2}} \ \mathbb{P}^3(\dj)
\end{equation}

In sections \ref{init}--\ref{def}--\ref{orb}--\ref{cat}, we essentially develop and discuss a number of aspects of the category of projective lines both in the finite and infinite case, and we show that the infinite case is quite different (through the construction of counter examples to theorems which hold true in the finite case). In section \ref{ext} we discuss extensions and embeddings (of lines into lines). This motivates us to consider {\em geometric closures} in section \ref{clos}. In the final section, which might well be the section with the most applications in its own right, we introduce and discuss AB-sets in order to construct rigid lines and rigid planes. Such geometries block the automorphisms of the ambient geometry in a precise way, and can be studied outside the category of axiomatic projective planes. The paper is complemented by a number of natural problems and questions which arise, and this aspect is one of our main motivations.

\section{Some initial definitions} 
\label{init}


A point-line incidence structure $\Upgamma = (\mP,\mL,\I)$ (where $\mP$ and $\mL$ are the point set and line set, and $\I$. is a symmetric incidence relation) is a {\em projective plane} if it satisfies the following axioms:
\begin{itemize}
\item[(a)]
any two distinct points are incident with one unique line; 
\item[(b)]
any two distinct lines are incident with one unique point;
\item[(c)]
there exist four distint points of which no three are incident with the same line.
\end{itemize}

Each projective plane has an {\em order} $N$; any of its points are then incident with $N + 1$ lines and any line is incident with $N + 1$ points.

A point-line incidence structure $\Upgamma = (\mP,\mL,\I)$ (where $\mP$ and $\mL$ are the point set and line set, and $\I$. is a symmetric incidence relation) is an  {\em affine plane} if it satisfies the following axioms:
\begin{itemize}
\item[(a)]
any two distinct points are incident with one unique line; 
\item[(b)]
if $A$ is a line and $a$ is a point not incident with $A$, then $a$ is incident with a unique line $B$ which does not meet $A$; 
\item[(c)]
there exist three distint points of which no three are incident with the same line.
\end{itemize}

Automorphisms and isomorphisms of and between projective and affine planes are defined in the standard way \cite{HP}. By $\Aut(\Upgamma)$, we denote the automorphism group of a  geometry $\Upgamma$. 

If $\Upgamma$ is a projective plane and $U$ is any line of $\Upgamma$, then by deleting the line $U$ and the points incident with $U$, one obtains an affine plane $\Upgamma_U$. On the other hand, if $\mA$ is an affine plane, then there exists a projective plane $\mP$ and a line $V$ in $\mP$ such that $\mP_V$ is isomorphic to $\mA$. The couple $(\mP,V)$ is unique up to isomorphism with respect to this property. Any automorphism of $\mP$ which fixes $V$ induces an automorphism of $\mP_V$, and every automorphism of $\mA$ defines an automorphism of $\mP$ which fixes $V$. If $\mA$ is an affine plane, then by $\overline{\mA}$ we denote its projective completion.

\section{Synthetic projective lines}
\label{def}

We are guided by the analogy between vector spaces (over division rings) and axiomatic affine planes with a distinguished point (which we compare to the zero vector). So consider any axiomatic affine plane $\mA$, and fix one arbitrary point $\omega$. Based on the discussion in section \ref{init}, we define the projective line $(\mA,\omega)$ as the set of lines $\mL_\omega$ of $\mA$ incident with $\omega$.  Note that we do not consider $\mL_\omega$ merely as a set: it comes with an embedding in $\mA$. The definition of the automorphism group of $(\mA,\omega)$ clearly illustrates this viewpoint. 

\subsection{Automorphism group}

An {\em automorphism} of $(\mA,\omega)$ is an element $\alpha$ of $\Aut(\mA)_\omega$, taking into account that $\alpha = \alpha'$ in $\Aut(\mA,\omega)$  with $\alpha, \alpha' \in \Aut(\mA)_\omega$ if and only if they induce the same action on $\mL_\omega$. 
It follows that  
\begin{equation}
\Aut(\mA,\omega)\ \cong\ \Aut(\mA)_\omega/N
\end{equation}
with $N$ the kernel of the action of $\Aut(\mA)_\omega$ on $\mL_\omega$|that is, the group of homologies with center $\omega$ and axis the line at infinity of $\mA$.

The notion of isomorphism is now obvious: we say that two synthetic lines $(\mA,\omega)$ and $(\mA,\omega')$ are {\em isomorphic} if there exists an isomorphism of planes 
\begin{equation}
\gamma:\ A\ \mapsto A'
\end{equation}
such that $\gamma(\omega) = \omega'$.

 In section \ref{cat} we will formally define the category of projective lines.

\subsection{Projective planes versus affine planes}

One could also consider the variation of the definition of synthetic projective line, in which the affine plane $\mA$ is replaced by a projective plane $\mP$. In that case, we could also dualize the definition, and see a projective type line as a line of a projective plane, which comes with an automorphism group which is induced by the group of the plane. 

We prefer the affine approach, on the one hand because of the vector space analogy, and on the other hand because it depends on the location of the base point (relative to the line at infinity, which is a natural model for the line).  Also, a number of very interesting problems arise in the affine approach, and these are among our motivations to study the problem. In fact, the very problems which arise when comparing the affine approach to the projective approach already are very promising, as we will see in this section. Let us remark that if $(\mA,\mu)$ is a projective line ($\mA$ an affine plane), then the embeddings 
\begin{equation}
(\mA,\mu)\ \hookrightarrow\ \mA\ \hookrightarrow\ \overline{\mA} 
\end{equation} 
implies that $\Aut(\mA,\mu)$ is a subgroup of $\Aut(\overline{\mA},\mu)$. 

The main two questions we ask are the following. Let $\mP$ be a projective plane, and consider, in the same way as above, the synthetic line (projective version) $(\mP,\mu)$ and $\Aut(\mP)_\mu/N$, with $N$ the kernel of the action of $\Aut(\mP)_\mu$ on the lines incident with $\mu$. For each line $U$ in $\mP$ which is not incident with $\mu$ we can consider the synthetic line $(\mP_U,\mu)$, with $\mP_U$ the affine plane arising from deleting $U$.  

\begin{question}[Inheritance]
\label{inher}
Is $\Aut(\mP)_\mu/N = \bigcup_{U \nI \mu}\Big({\Aut(\mP)}_{U,\mu}\Big/M_U\Big)$, with $M_U$ the kernel of the action of  ${\Aut(\mP)}_{U,\mu}$ on the lines incident with $\mu$. 
\end{question}

In other words: can we recover each automorphism of the projective type line $(\mP,\mu)$ from some automorphism of an affine type line $(\mP_U,\mu)$?\\

We will construct counter examples in the infinite case, which will indicate that many counter examples exist. In the finite case, many counter examples exist as well, and we will explain some thoughts on that matter. A second, highly related but altogether different question which we need to understand is the following:  

\begin{question}[Generation]
Is $\Aut(\mP)_\mu/N = \Big\langle \Big({\Aut(\mP)}_{U,\mu}\Big/M_U\Big)\ |\ U \nI \mu \Big\rangle$, with $M_U$ the kernel of the action of  ${\Aut(\mP)}_{U,\mu}$ on the lines incident with $\mu$. 
\end{question}

Here we ask whether each automorphism of the projective type line $(\mP,\mu)$ is {\em generated} by automorphisms of affine type lines $(\mP_U,\mu)$.

\subsection{Some general remarks}

We first recall the next theorem:

\begin{theorem}[Hughes and Piper \cite{HP}, Theorem 13.3]
\label{fix}
An automorphism of a finite projective plane has an equal number of fixed points and lines.
\end{theorem}

Now any automorphism $\alpha \in \Aut(\mP)$ which fixes $\mu$, also fixes some line $M$. If $M$ is not incident with $\mu$, then $\alpha$, as an automorphism of $(\mP,\mu)$, is inherited from an automorphism of $(\mP_M,\mu)$. But theoretically, automorphisms which only fix a flag of $\mP$, or more generally which fix $m < 0$ lines on a point $\mU$ and $n > 0$ points on exactly one these lines (and no other lines and points), could exist, and would not be inherited. (Call such automorphisms ``of projective origin.'') We will now construct counter examples which work both in the finite and infinite case for Desarguesian planes, and more generally for semifield planes. In the next subsection, we will construct a rather different class of counter examples. 

Call a projective plane $\mP$ a {\em semifield plane} if there exists a flag $(x,Y)$ for which $x$ is a translation point (with translation group $T(x)$), and $Y$ a translation line (with translation group $T(Y)$). Now let $t \in T(Y)^\times$ be such that it fixes the point $u \ne x$ on $Y$ linewise (an element of $T(Y)$ which is not an $(x,Y)$-elation has this property). Let $\widetilde{y}$ be an element in $T(x)^\times$ which does not fix $Y$ pointwise (again, any element which is not an $(x,Y)$-elation is such an element). Note that $T(Y) \cap T(x)$ is precisely the complete group of all $(x,Y)$-elations. Then obviously $t\widetilde{t}$ fixes the flag $(x,Y)$ and no other points and lines.  

For Desarguesian projective planes, the converse is also true in case of linear automorphisms. For, let $(x,Y)$ be a flag in such a plane $\mP$, and let $\alpha$ be a linear automorphism fixing $x$, $Y$ and no other points and lines. Suppose $\alpha$ induces the (linear) automorphism $\alpha_Y$ on $Y$; then it is easy to see that $T(x)$ contains an element 
$\widehat{\alpha}$ which induces $\alpha_1$. In the same way, if $\alpha_x$ is the linear automorphism induced by $\alpha$ on the lines incident with $x$, 
we take an element $\widetilde{\alpha}$ in $T(Y)$ which induces $\alpha_x$ on these lines. Now consider $\epsilon := \alpha(\widehat{\alpha}\widetilde{\alpha})^{-1}$; this is 
an automorphism which fixes $x$ linewise and $Y$ pointwise, so $\epsilon \in T(x) \cap T(Y) \leq \langle T(x), T(Y) \rangle = T(x)T(Y)$. It follows that $\alpha$ is also contained 
in $T(x)T(Y)$. If $\alpha$ is not linear, this might depend on the field over which the plane is defined. Suppose that $\mP$ is finite and defined over $\F_q$, with $q = p^h$ and 
$p$ a prime. Suppose that $(p,h) = 1$, so that $p$ does not divide $\vert \texttt{Gal}(\F_q/\F_p) \vert$. As the order of $\alpha$ must be a power of $p$ (otherwise we 
have more fixed points and fixed lines), and since $[ \mathsf{P\Upgamma L}_3(q) : \mathsf{PGL}_3(q) ] = h = \vert \mathsf{Gal}(\F_q/\F_p) \vert$, it follows that $\alpha$ is linear.        





\subsection{Infinite case}

In the infinite case, we need a different approach. For, suppose that $\gamma$ is an automorphism of the infinite projective plane $\mP$.
Let $\mP = \{ x_i\ \vert\ i \in \mI \}$ be the set of points of the plane, and $\mL = \{ M_i \ \vert\ i \in \mI \}$ be the set of lines. Let $A$ be an incidence matrix  for $\mP$ ($a_{ij} = 1$ if and only if $x_i$ is incident with $M_j$, otherwise it is $0$). 
Define a matrix $P = (p_{ij})$ as follows: $p_{ij} = 1$ if $x_i^\gamma = x_j$; (otherwise $0$) and define a matrix $L = (l_{ij})$ in a similar way: $l_{ij} = 1$ if $l_i^\gamma = l$ (otherwise $0$). Note that $P$ and $L$ are infinite permutation matrices. Now consider $(PA)_{ij} = \sum_{k \in \mI}p_{ik}a_{kj} =: a_{uj}$, where $u$ is uniquely determined by the fact that $p_{iu} = 1$: $x_i^\gamma = x_u$. Next, we look at the $(i,j)$-entry $(AL)_{ij}$ of $AL$, which is given by  $\sum_{k \in \mI}a_{ik}l_{kj} = a_{iv}$, where $v$ is uniquely determined by the fact that $M_v^{\gamma} = M_j$. Now $a_{uj} = 1$ if and only if $x_u \I M_j$ if and only if $x_u^{\gamma^{-1}} \I M_j^{\gamma^{-1}}$ if and only if $x_i \I M_v$ if and only if $a_{iv} = 1$, implying that 
\begin{equation}
PA\ =\ AL.
\end{equation}

In the finite case one now easily concludes that $P$ and $L$ have the same trace, translated into the fact $\gamma$ has the same number of fixed points as fixed lines. But in the infinite case this is essentially where the story ends. In the next subsection, we construct counter examples to Theorem \ref{fix} in the infinite case.

\subsection*{Counter examples} 

Consider the Fano plane $\Pi$, and let $(x,M)$ be an anti-flag. Let $\alpha$ be an automorphism of $\Pi$ of order $3$ which fixes $x$, $M$ and no other point or line, so that point-orbits and line-orbits of $\langle \alpha \rangle$  which do not contain $x$ nor $L$ all have size $3$ (and note that such automorphisms exist). Now remove $M$ and the points incident with $M$ (which form one point-orbit) to obtain a point-line configuration $\mC$  (affine plane of order $2$), and $\alpha$ induces an automorphism of $\mC$ which we denote by $\alpha_0$. 

We now perform a free construction on $\mC$ to obtain a countably infinite free projective plane $\overline{\mC}$.
Put $\mC = \mC_0$. In step one, we add new points to $\mC_0$ for each pair of nonconcurrent lines in $\mC_0$ to obtain $\mC_1$. In the next step, we add new lines for each pair of noncollinear points to obtain $\mC_2$.  And so on. The free plane $\overline{\mC}$ is 
\begin{equation}
\cup_{i \in \mathbb{N}}\mC_i.
\end{equation}  
In each step, we also construct an automorphism of the newly obtained point-line configuration. This is easy. For instance, one we have constructed $\mC_1$, the new points naturally form an orbit of size $3$ (by looking at the image of the lines) and we obtain $\alpha_1$. The same is true for $\mC_2$, and we obtain $\alpha_2$. By induction, we obtain that each $\alpha_i$ is well defined, and its ``union" $\overline{\alpha}$ is an automorphism of $\overline{\mC}$ which fixes the point $x$, and all other points and lines live in orbits of $\langle \overline{\alpha} \rangle$ of size $3$. (Note that if the new line $N$ is added in stage $\mC_m$, then it only contains two points which both live in separate orbits of size $3$, so that $N$ cannot be fixed, and it lives in an $\langle \alpha \rangle$-orbit of size $3$.)

So there are no fixed lines.  \\

Note that the automorphism $\overline{\alpha}$ cannot be inherited by some affine type projective line automorphism based at $x$: 

\begin{theorem}
\label{inhercount}
The automorphism $\overline{\alpha}$ defines a nontrivial automorphism of the projective type line $(\overline{\mC},x)$ which is not inherited by any affine type line automorphism. 
\end{theorem}

{\em Proof}.\quad 
Suppose there is some line $U$ in $\overline{\mC}$ not incident with $x$, so that there is an automorphism of $\overline{\mC}_U$ (which we see as an automorphism of $\overline{\mC}$ which fixes $U$) which induces $\overline{\alpha}$ on the lines of $\overline{\mC}$ incident with $x$. Then $\beta := \overline{\alpha}\gamma^{-1}$ is a central collineation of $\overline{\mC}$ with center $x$. Now let $\Delta$ be a triangle with sides $A, B, C$, such that $x$ is not incident with $A, B$ and $C$. Then $A, B, C, A^\beta, B^\beta, C^\beta$ define a Desargues configuration $\mD$. Now $\mD$ is a finite point-line configuration with three points per line and three lines per point, so by definition the configuration is {\em confined}. By Heartshorne \cite{HH}, it follows that $\mD^{\overline{\alpha}}$ is a subgeometry of the begin configuration $\mC_0$, which is the contradiction we were seeking. So $\overline{\alpha}$ indeed defines a nontrivial automorphism of the projective type line $(\overline{\mC},x)$ which is not inherited by any affine type line automorphism in $\overline{\mC}$. \eop \\

Applied to Question \ref{inher}, Theorem \ref{inhercount} yields counter examples in the infinite case. \\

Of course, this method immediately generalizes to the following theorem (using the same notation). 

\begin{theorem}
Let $\Pi$ be a finite projective plane, and let $(x,M)$ be an anti-flag. Suppose $\alpha$ is an automorphism of $\Pi$ which fixes $x$, $M$ and no other points nor lines. Remove $M$ and its points to obtain the affine plane $\mC$, and free construct the projective plane $\overline{\mC}$. Then $\overline{\alpha}$ is an automorphism of $\overline{\mC}$ which fixes $x$, and no other points and lines of $\overline{\Pi}$. \eop 
\end{theorem}

\subsection{Group theoretical discussion}

On the abstract group-theoretical level, we pose the following related problem. 

\begin{question}
Let $(H,X)$ be a permutation group. When is $H$ generated by (all) its subgroups $H_x$, where $H_x$ is the stabilizer of $x \in X$. 
\end{question}

In \cite{MO} it is shown that if $H$ is transitive on $X$, then it will be generated by its point-stabilizers if and only if it does not have a proper system of imprimitivity upon which it acts sharply transitively. The proof is very simple: let $N$ be the normal subgroup of $H$ generated by all the point-stabilizers $H_x$. If $N \ne H$, then clearly $N$ cannot act transitively on $X$ (as it contains point-stabilizers). If $\Omega$ is an $N$-orbit, then $H_\Omega$ and $N$ act transitively on $\Omega$, so that $H_\Omega = NH_x$ for some $x \in \Omega$. So $H_\Omega = N$, and $H_\Omega$ fixes all $N$-orbits. Conversely, let $\mB$ be a proper system of imprimitivity upon which $H$ acts sharply transitively, and let $K$ ($\ne H$) be the kernel of this action. Let $x \in X$, and let $\Omega$ be the block in $\mB$ containing $x$. Then $H_x \leq H_\Omega = K$. As $K$ is a normal subgroup of $H$, and as all point-stabilizers are conjugate to $H_x$, it follows that $K$ is the subgroup of $H$ generated by all point-stabilizers.

\section{Orbits and lines}
\label{orb} 

Many nonisomorphic synthetic projective lines can arise from one and the same plane $\mP$, and the isomorphism classes are obviously determined precisely by the orbits of $\Aut(\mP)$ on the set $\mho(\mP)$ of anti-flags in $\mP$. Rigid planes (see also section \ref{rig}), for example, yield the maximal number of nonisomorphic projective lines. 
Here is the extremal case on the other end of the spectrum.

\begin{theorem}
Let $\mP$ be a finite  projective plane. Then all lines $(\mP_U,\omega)$ (with $U$ a line in $\mP$ and $\omega$ not incident with $\omega$) are isomorphic 
if and only if $\mP$ is Desarguesian.
\end{theorem}

{\em Proof}.\quad
Obviously this is the case if and only if $\mP$ acts transitively on anti-flags. In that case, it follows from Ostrom and Wagner \cite{OstWag} (see also Cameron and Kantor \cite{CamKan}).
\eop \\ 

In the infinite case, no such conclusion can hold due to model-theoretic considerations. In \cite{Tent} Tent shows the following theorem (which we only state for projective planes, but which she obtains for all generalized polygons). We refer to her paper for the notational use. 

\begin{theorem}[Tent \cite{Tent}]
\label{tent}
There exists a countable projective plane $\Upgamma$ which is not isomorphic to $\mathbb{P}^2(\dj)$ for any division ring $\dj$, which contains all finite (known and unknown) projective planes (up to isomorphism).  
Its automorphism group $\Aut(\Upgamma)$ acts transitively on all finitely generated $\mL_3$-structures (in $\Upgamma$) of given isomorphism type. In particular, 
\begin{itemize}
\item[(a)]
$\Aut(\Upgamma)$ is a group with a BN-pair; 
\item[(b)]
$\Aut(\Upgamma)$ acts highly transitively on (all lines incident with) every point, and on (all point incident with) every line;
\item[(c)]
if $\upgamma = (x_0,x_1,\ldots,x_6 = x_0)$ is a $6$-cycle in the incidence graph of $\Upgamma$, then the elementwise stabilizer of $\upgamma$ acts highly transitively on $\Upgamma_1(x_1) \setminus \{ x_0, x_2 \}$ (and obviously $x_1$ is arbitrary in this context). 
\end{itemize}
\end{theorem}

Here, $\Upgamma_1(x_1)$ denotes the set of vertices adjacent to (and different from) $x_1$ (in the incidence graph).  

The fact that $\Aut(\Upgamma)$ is a group with a BN-pair translates into the fact that its acts transitively on the ordered triangles in $\Upgamma$. As a by-product, it acts 
transitively on the anti-flags, so that $\Upgamma$ is a nonclassical projective plane for which all projective lines $(\Upgamma^U,\omega)$ (where $(\omega,U)$ is any anti-flag in $\Upgamma$) are isomorphic, where $\mA = \Upgamma^U$ is the affine plane obtained by deleting the line $U$. 

Projective planes $\mP$ in which all projective-type projective lines $(\mP,\omega)$ are isomorphic are precisely those planes with a point-transitive automorphism group. In the finite case, there are several long-standing conjectures which imply that such planes are Desarguesian. In the infinite case such results cannot be true, as we have seen in Theorem \ref{tent}.   \\

Note that if we fix some line $U$ in a projective plane $\mP$, the isomorphism class of $(\mP,\omega)$, with $\omega$ a point not on $U$ solely depends on the point-orbits of $\Aut(\mP)_U$. So one line in the plane $\mP$ can define many nonisomorphic affine-type projective lines. In general, this relative viewpoint is very different than in the case of 
``classical'' projective lines.

\section{Categorical setting}
\label{cat}

Define the category $\wp$ as follows. Its object class $\ob(\wp)$ consists of synthetic projective lines $(\mA,\omega)$, and a morphism 
\begin{equation}
f:\ (\mA, \omega)\ \mapsto\ (\mA', \omega')
\end{equation}
is simply a morphism of affine planes (points are mapped to points, lines to lines, incidence is preserved, parallel classes are sent to parallel classes) for which $f(\omega) = \omega'$.

\subsection{Compatible automorphisms}

Note that there is no hope that automorphisms $\upgamma$ of $(\mA,\omega)$ survive the following diagram:  

\begin{center}
\begin{tikzpicture}[>=angle 90,scale=2.2,text height=1.5ex, text depth=0.25ex]
\node (a0) at (0,3) {$(\mA,\omega)$};
\node (a1) [right=of a0] {$(\mA',\omega')$};

\node (b0) [below=of a0] {$\Aut(\mA,\omega)$};
\node (b1) [below=of a1] {$\Aut(\mA',\omega')$};

\draw[->,font=\scriptsize,thick]
(a0) edge node[left] {} (b0)
(a1) edge node[right] {} (b1);

\draw[->,font=\scriptsize,thick]
(a0) edge node[auto] {$\upgamma$} (a1)
(b0) edge node[auto] {$\widetilde{\upgamma}$} (b1);


\end{tikzpicture}
\end{center}

Only automorphisms $\alpha$ of $\mA$ ``commuting with $\upgamma$'' have an interpretation as an automorphism on $\upgamma(\mA)$: if $X, Y \in \upgamma^{-1}(Z)$ for each $Z$, point or line in $\upgamma(\mA)$, then $X^\alpha, Y^\alpha \in \upgamma^{-1}(\widehat{Z})$, where $\widehat{Z} = X^{\alpha\upgamma}$. Call such automorphisms $\alpha$ {\em compatible} if both $\alpha$ and $\alpha^{-1}$ have this property. (When the order of $\alpha$ is finite, as we will see below, $\alpha^{-1}$ automatically has this property if $\alpha$ does.) 
But is easy to construct ``noncompatible automorphisms,'' as we will show now. We first recall a result of Pasini \cite{P} on epimorphisms of generalized polygons, which generalizes an older result of Skornjakov \cite{Skorn} and Hughes \cite{H} (the latter two handling the case of projective planes).

\begin{theorem}[Skornjakov--Hughes--Pasini \cite{Skorn,H,P}]
\label{HP}
Let $\alpha$ be a morphism from a thick (possibly infinite) generalized $m$-gon $\mE$ to a thick (possibly infinite) generalized $m$-gon $\mE'$, with $m \geq 3$. If $\alpha$ is surjective, then 
either $\alpha$ is an isomorphism, or each element in $\mE'$ has an infinite fiber in $\mE$.
\end{theorem}

Thas and Thas handle the ``thin case'' in \cite{JATKTepi}:

\begin{theorem}[Thas and Thas \cite{JATKTepi}]
\label{GT}
Let $\Phi$ be an epimorphism of a thick projective plane $\mP$ onto a thin projective plane $\Delta$ of order $(1,1)$.  Then exactly two classes of epimorphisms $\Phi$ occur (up to a suitable permutation of the points of $\Delta$), and they are described as follows. 

\begin{itemize}
\item[{\rm(a)}] 
The points of $\Delta$ are $\overline{a}, \overline{b}, \overline{c}$, with $\overline{a} \sim \overline{b} \sim \overline{c} \sim \overline{a}$, and put $\Phi^{-1} (\overline{x}) = \wt X$, with $\overline{x} \in \lbrace \overline{a}, \overline{b}, \overline{c} \rbrace$.

   Let $(\wt{A}, \wt{B})$, with $\wt{A} \ne \emptyset \ne \wt{B}$, be a partition of the set of all points incident with a line $L$ of $\mP$. Let  $\wt{C}$ consist of the points not incident with $L$. Furthermore, $\Phi^{-1}(\overline{a}\overline{b}) = L$, $\Phi^{-1}(\overline{b}\overline{c})$ is the set of all lines distinct from $L$ but incident with a point of $\wt{B}$ and $\Phi^{-1}(\overline{a}\overline{c})$ is the set of all lines distinct from $L$ but incident with a point of $\wt{A}$.
   \item[{\rm(b)}] 
   The dual of (a).
\end{itemize}
\end{theorem}

Consider any thick projective plane $\mP$, and let $\alpha$ be a nontrivial automorphism of $\mP$ (this is an assumption). Then there is some line $U$ which is not fixed by $\alpha$. 
Let $a$ and $b$ be different points of $\mP$ on $U$, and suppose that $a^\alpha \not\I U$ while $b^\alpha \I U$ (note that such couples $(a,b)$ exist: just let $b$ be the point which maps onto $U \cap U^\alpha$).  
Now let $\upgamma: \mP \mapsto \mP'$ be a morphism of projective planes, where $\upgamma(\mP) =: \Delta$ is any triangle in $\mP'$. Construct $\upgamma$ in such a way that 
 $\{ U \} = \upgamma^{-1}(V)$, where $V$ is a side of $\Delta$, and such that $a$ and $b$ are contained in the same point-fiber of $\upgamma$. Then $a^\alpha$ and $b^\alpha$ are not contained in the same point-fiber of $\upgamma$, so $\alpha$ and $\upgamma$ are not compatible (once one descends to the appropriate affine planes). 

\begin{observation}
Let $C(\upgamma) \subseteq \Aut(\mP)$ be such that each element of $C(\upgamma)$ is compatible with the given morphism of planes $\upgamma: \mP \mapsto \mP'$. Then  $C(\upgamma)$ is a subgroup of $\Aut(\mP)$, and $\upgamma^{-1}C(\upgamma)\upgamma \cong C(\upgamma)$ is an automorphism group of $\mP'$.  
\end{observation}

{\em Proof}.\quad
Let both $\alpha$ and $\beta$ be compatible with $\upgamma$, and take $X, Y \in \upgamma^{-1}(Z)$ for some element $Z = X^\upgamma$ in $\upgamma(\mP)$. Then $X^\alpha, Y^\alpha \in \upgamma^{-1}(X^{\alpha\upgamma})$, so ${(X^{\alpha})}^\beta, {(Y^{\alpha})}^\beta \in \upgamma^{-1}\Big({(X^{\alpha})}^{\beta\gamma}\Big)$. As $\alpha^{-1} \in C(\upgamma)$ for each $\alpha \in C(\upgamma)$, it follows that $C(\upgamma)$ is a subgroup of $\Aut(\mP)$. The second part is straightforward (note that if $X \in \upgamma^{-1}(Z)$, then $X^{\alpha\gamma}$ is independent of the choice of $X$ in the fiber, so elements of $\upgamma^{-1}C(\upgamma)\upgamma \cong C(\upgamma)$ are well defined). 
\eop

\subsection{Morphisms and sublines}

Note that if $\gamma:\ (\mA, \omega)\ \mapsto\ (\mA', \omega')$ is a morphism, then $\gamma(\mA,\omega)$ is a (possibly degenerate) ``subline'' of $(\mA',\omega')$. This follows from the fact that $\gamma(\overline{\mA})$ (where $\overline{\mA}$ is a projective closure of $\mA$ and $\gamma(\overline{\mA})$ its natural image) is a (possibly degenerate) projective plane (if $U, V$ are distinct lines in $\gamma(\overline{\mA})$, then if $U' \in \gamma^{-1}(U)$ and $V' \in \gamma^{-1}(V)$, the image 
of $U' \cap V'$ under $\gamma$ is an intersection point of $U$ and $V$).

\section{Extensions and embeddings}
\label{ext}

Let $(A,\omega)$ be a synthetic projective line. We say that $(A',\omega')$ is an {\em extension} of $(A,\omega)$, or that $(A,\omega)$ is a {\em subline} 
of $(A',\omega')$ if there is an injective morphism 
\begin{equation}
\rho:\ A\ \mapsto A'
\end{equation}
such that $\rho(\omega) = \omega'$. Since $A \cong \rho(A)$, we can see $A$ as a subplane of $A'$. Note that there is no reason to think that $\Aut(\mA,\omega)$ is a subgroup 
of $\Aut(\mA',\omega')$. In general, we expect these groups to be unrelated, although, for example, in case $\mA$ is a Desarguesian plane over $\F_q$ and $\mA'$ a 
Desarguesian plane over $\F_{q^n}$, $n \geq 1$, we indeed have a natural inclusion of subgroups.  

\subsection{Desarguesian sublines}

In the classical context, the notion of subline coincides with a line over a subfield.

\begin{theorem}
\label{des}
Let $(A,\omega)$ be a subline of $(A',\omega')$, where $A' = \A^2(\F_q)$ is a Desarguesian affine plane over $\F_q$. Then $U$ (line of infinity of $\rho(A)$) is a line 
defined over a subfield of $\F_q$, that is, a projective of $U'$ (line of infinity of $A'$) in the classical sense. 
\end{theorem}

{\em Proof}.\quad
Let 
$\rho:\ A\ \mapsto A'$
be an injective morphism such that $\rho(\omega) = \omega'$. We have that $\rho(A)$ is a subplane of $A' \cong \A^2(\F_q)$, so $\rho(A)^U$ is a projective subplane of 
${A'}^{U'} \cong \bP^2(\F_q)$. It follows immediately that $\rho(A)^U \cong \bP^2(\ell)$ with $\ell$ some subfield of $\F_q$, and the statement follows. \eop \\


\subsection{Self-embeddings} 

Note that in the infinite case, we can easily ``self-embed'' projective lines. Here is an example. 

Consider the classical plane $\bP^2(\ell)$, $\ell$ a field. We work with homogeneous coordinates $(x : y : z)$. Let $U$ be the line $z = 0$, and let $\omega$ be the point $(0, 0, 1)$. Note that each isomorphism of fields $\alpha: \ell \cong \ell'$ gives rise to an isomorphism of planes sending $(0,0,1)$ to $(0,0,1)$ and sending $z = 0$ to $z = 0$; ``the'' isomorphism is given by 
\begin{equation}
(a, b, c) \mapsto (a^\alpha, b^\alpha, c^\alpha).   
\end{equation} 

Now put $\ell = k(x_1,x_2,x_3,\ldots)$, and define $\alpha$ by acting as the identity on $k$, and sending
\begin{equation}
x_1\ \mapsto\ x_2\ \mapsto\ x_3\ \mapsto \cdots 
\end{equation} 
Then $\ell = k(x_1,x_2,\ldots)$ is sent to $\ell' = k(x_2,x_3,\ldots)$, giving rise to a strict embedding of $(\bP^2(\ell')_U,\omega) \cong (\bP^2(\ell)_U,\omega)$ in $(\bP^2(\ell)_U,\omega)$. 

\subsection{Remark} 

More generally, let $\gamma$ be an isomorphism of $L = (\mA,\omega)$ onto $L' = (\mA',\omega')$ while $L' \ne L$ is a subline of $L$. Using the sequence 
\begin{equation}
\cdots \ \mapsto\ L^{\gamma^{-2}} \ \mapsto\ L^{\gamma^{-1}} \ \mapsto\ L \ \mapsto\ L^{\gamma}  \ \mapsto\ L^{\gamma^{2}} \ \mapsto \cdots
\end{equation}
one can now  construct larger lines through the union / direct limit of this system. We want to explore this idea a little further in the next section.

\subsection{Field analogy}

Considering Theorem \ref{des}, there obviously is an interesting analogy between line extensions and field extensions. Note that any line $(\mA,\mu)$ is embedded in (many) other lines: just add an extra point $x$ to $\overline{\mA}$, and construct the free closure of $\overline{\mA} \cup \{ x\}$ to obtain a projective plane $\Upgamma$. If $U$ is the line at infinity of $\mA$, then $(\Upgamma_U,\mu)$ is an extension of $(\mA,\mu)$. In the next section, we will make some further remarks on the aforementioned analogy.

\section{Geometric closure}
\label{clos}

Call a line $(\mA,\omega)$ {\em finite} if $\mA$ is finite. Call an affine plane $\mA$ {\em set-like} if its collections of points and lines are sets; otherwise 
it is {\em class-like}. We use the same terminology for lines.

Let $(\mA,\omega)$ be a set-like synthetic projective line, and let $C(\mA,\omega)$ be the category of all set-like synthetic projective lines which contain 
$(\mA,\omega)$. (Note that all the planes have the same line at infinity.) 
In a natural way we can make $C(\mA,\omega)$ into a direct system (by taking natural embeddings as connection morphisms). If $(\mB,\omega)$ and $(\mC,\omega)$ are elements of $C(\mA,\omega)$, then the free completion of $\mB^U \cup \mC^{U}$ (where $U$ is the common line at infinity) defines an upper bound of both lines after removing $U$. 
We call $\overline{(\mA,\omega)} := \varinjlim{C(\mA,\omega)}$ the {\em geometric closure} of $(\mA,\omega)$. Note that $\overline{(\mA,\omega)}$ is a class-like line, and 
that $\overline{\mA}$ (the obtained affine plane) also is class-like. 

Note that if $(\mB,\omega) \in C(\mA,\omega)$, we have that 
\begin{equation}
\overline{(\mA,\omega)}\ = \ \overline{(\,B,\omega)},
\end{equation}
showing a striking similarity with the notion of algebraic closure. Also note that $\overline{(\mA,\omega)}$ has the same line at infinity as $(\mA,\omega)$. 

However, geometric closures of lines of type $\bP^2(k)$ are obviously much ``larger'' than lines of type $\bP^2(\overline{k})$, where $\overline{k}$ is an algebraic closure of $k$. We could modify the construction by restricting $C(\mA,\omega)$ as follows:
\begin{itemize}
\item
If $(\mA,\omega)$ is finite, only consider finite and countable extensions of $(\mA,\omega)$.
\item
If $(\mA,\omega)$ is not finite, we only consider extensions $(\mB,\omega)$ such that $\mA$ and $\mB$ have the same order. 
\end{itemize} 

Even in the modified variation, if we would start from $(\A^2(\F_q),\omega)$, then 
\begin{equation}
\overline{(\A^2(\F_q),\omega)}\ \not\cong \ (\A^2(\overline{\F_q}),\omega).
\end{equation}

This is easy to see: $C(\A^2(\F_q),\omega)$ contains nonclassical planes; just freely complete $\A^2(\F_q)^U$ and an extra point, and remove $U$ again: we obtain an element in $C(\A^2(\F_q),\omega)$ which is not a subline of $(\A^2(\overline{\F_q}),\omega)$. On the other hand, $\overline{(\A^2(\F_q),\omega)}$ contains $(\A^2(\overline{\F_q}),\omega)$ as a subline. 

\medskip
\bq
Can one geometrically define a notion of geometric closure which coincides with the notion of algebraic closure for finite Desarguesian lines?
\eq

The following variation is equally interesting. 

\medskip
\bq
Can one geometrically define a notion of geometric closure which coincides with the notion of algebraic closure for finite Desarguesian planes?
\eq

Here the term ``geometrically'' must be interpreted at a synthetic-geometrical level (and not algebro-geometric). One way to approach the problem is to consider $\bP^2(\overline{\F_q})$ for the finite field $\F_q$, and use the fact that $\overline{\F_q}$ is the union of its finite subfields, noting that for each positive integer $n \ne 0$, $\overline{\F_q}$ contains one 
unique subfield isomorphic to $\F_{p^n}$, with $q$ a power of the prime $p$.


\medskip
\section{Rigid projective lines, realization and AB-sets}
\label{rig}

The following question is natural. 

\begin{question}
\label{qgroup} 
Let $G$ be a given abstract group. Does there exist a projective or affine type projective line for which $G$ is isomorphic to the full automorphism group? 
\end{question}

Typically we are interested in (multiply) transitive groups (think of groups with a BN-pair of rank $1$ or even general $2$-transitive groups), but also in the other side 
of the spectrum: the trivial group. Before passing to the latter case, we will make some remarks about some specific transitive examples. We first introduce a specification of 
Question \ref{qgroup}. Let $\bP$ be a projective line in any of the other theories (vector spacial, split BN-pairs of rank $1$, algebro-geometric / topological, etc.), and 
suppose $\Aut(\bP)$ is its automorphism group in that theory. Can we find a projective-type or affine-type projective line which defines equivalent permutation groups? We call 
lines such $\bP$ {\em realizable}.

\subsection*{Sharply $2$-transitive groups}

In the finite case, synthetic lines with a sharply $2$-transitive automorphism group can hardly be realized upon accepting the prime 
power conjecture. We first recall the next theorem from Zassenhaus (cf. \cite[Theorem 20.3]{Pass}).  

\begin{theorem}
If $(H,X)$ is a finite sharply $2$-transitive permutation group, then $\vert X \vert$ is a prime power $q$, and either $H$ acts as a one-point-stabilizer in $\mathbf{P\Upgamma L}_2(q)$ in its natural action on the projective line $\bP^1(\F_q)$, or $q \in \{  5^2, 7^2, 11^2, 23^2, 29^2  \}$.  
\end{theorem}

It follows that if $(H,X)$ is a finite (faithful) sharply $2$-transitive permutation group, then $\vert X \vert$ is a prime power, so if $(\mA,\omega)$ is a finite projective line with a sharply $2$-transitive automorphism group, where $\mA$ has order $N$, then $N + 1 = p^m$ for some prime $p$ and positive integer $m$. If the prime power conjecture would be true, we end up with couples $(N,N + 1)$ where both entries are prime powers. It is well known that for such couples either $N$ or $N + 1$ are powers of $2$, and if $(N, N + 1) \ne (8,9)$, the other one is a prime number. In the infinite case there is much more elbow room. 

\subsection*{The general linear groups}

Obviously, the groups $\mathsf{P\Gamma L}_2(\dj)$, with $\dj$ a division ring, are realizable.

\subsection*{Highly transitive lines}

Let $\Upgamma$ be the plane constructed by Tent as in Theorem \ref{tent}. Then for any point $\omega$, we have that $\Aut(\Upgamma,\omega)$ acts highly transitively on the points of $(\Upgamma,\omega)$. Now let $U$ be any line of $\Upgamma$, and $v$ any  point not incident with $U$. As $\Aut(\Upgamma)$ is a group with a BN-pair, it follows that $\Aut(\Upgamma)_{(U,v)}$ acts $2$-transitively on the lines incident with $v$. Let $\mA$ be the affine plane obtained by deleting $U$ from $\Upgamma$. Then by Theorem \ref{tent}(c), it follows that $\Aut(\mA,v)$ acts highly transitively on the points of $(\mA,v)$. \\

We now pass to the other side of the spectrum. 

\subsection{Rigid lines}

Call a projective line $(\mA,\omega)$ {\em rigid} if its automorphism group is the trivial group. This immediately translates into the fact that $\Aut(\mA)_\omega$ is a subgroup of the group of homologies of $\mA$ with center $\omega$ (and axis the line at infinity). Although it is a wide-open question whether finite rigid projective planes exist, we will show that as a direct corollary of known results, finite (and infinite) rigid affine and projective type lines exist. Especially in the finite case, this is particularly interesting. 

\subsection{Kantor's translation planes}
\label{Kanpl}

In \cite{Kan}, Kantor constructed infinite classes of finite translation planes $(\pi,U,T)$ of which the full automorphism group $\Aut(\pi)$ coincides with the translation group $T$ (and hence have kernel isomorphic to $\F_2$). More precisely, we have such planes for each prime power $q = 2^h$ with $h$ odd, composite and $h \geq 3$. 
For each point $\omega$  in $\pi \setminus T$, $(\pi_U,\omega)$ defines a rigid projective line. The line $U$ itself defines a projective type projective line (it is the only line in $\pi$ with this property). 

The following problem might be interesting:

\begin{problem} 
Study affine planes $\mA$ (with line at infinity $U$) which yield a non-rigid projective type line $U$, but for every point $\omega$ a rigid projective line $(\mA,\omega)$. 
\end{problem}

This translates into the fact that $\Aut(\mA)$ is nontrivial and does not fix every parallel class, but if $\alpha \in \Aut(\mA)$ fixes an affine point, then all parallel classes {\em are} fixed.

\subsection{Rigid planes, symmetry breaking and AB-sets} 

Obviously it is much harder to produce rigid projective (and affine) planes than rigid lines. Note that a rigid projective plane yields a rigid affine plane after deleting any line; 
the projective completion of a rigid affine plane is not necessarily rigid though (the Kantor planes of the previous section are finite examples).  
In this subsection we have a closer look at an easy construction of infinite rigid projective planes through so-called ``AB-sets."\\

It is an open question as to whether finite rigid projective planes exist. In the infinite case, things are very different: Mendelsohn \cite{Mendel} showed that for every group $G$, there is a projective plane $\mP$ such that its full automorphism group is isomorphic to $G$. The case $G = \{ \id \}$ yields rigid projective planes, and hence for every anti-flag $(a,A)$ in such a plane $\mP$, we have that $(\mP_A,a)$, respectively $(\mP,A)$, is a rigid affine type, respectively projective type, projective line. 

It is not hard to construct infinite rigid projective planes, and so we will describe a procedure with which one can do this. 

\begin{itemize}
\item
Start with a point-line geometry $\Upgamma$ which does not violate the axioms of a projective plane. (In this context, it will be most convenient to have $\Upgamma$ finite.)     
\item
Leave away a point-line subconfiguration $\mC$ of $\Upgamma$ such that the new geometry $\Upgamma \setminus \mC =: \mC_0$ is confined and rigid. 
\item
Free complete $\mC_0$ to obtain a rigid projective plane $\overline{\mC_0}$.
\end{itemize}

Why is the plane $\overline{\mC_0}$ rigid? First note that if an automorphism of $\overline{\mC_0}$ fixes $\mC_0$ elementwise, then by its very definition, $\overline{\mC_0}$ is fixed elementwise by the completion of $\alpha$. Now let $\gamma$ be an automorphism of $\overline{\mC_0}$; then $\mC_0^\gamma$ is also confined, so by \cite{HH} we must have $\mC_0^\gamma \leq \mC_0$. As $\mC_0$ is finite by definition, it follows that $\gamma$ induces  an automorphism, hence the identity, on $\mC_0$.   

We propose to define {\em automorphism-blocking sets} (``AB-sets'') as follows. With $\Upgamma$ a point-line geometry, we say that $\mC$ is an {\em automorphism-blocking set} of $\Upgamma$ (or that $\mC$ {\em blocks the automorphisms} of $\Upgamma$, or that it {\em breaks the symmetry} of $\Upgamma$), if the following conditions are satisfied:
\begin{itemize}
\item[AB1]
every automorphism of $\Upgamma \setminus \mC$ extends to an automorphism of $\Upgamma$;
\item[AB2]
the elementwise stabilizer $\Aut(\Upgamma)_{[\mC]}$ of $\mC$ in $\Aut(\Upgamma)$ is trivial;
\item[AB3]
$\Aut(\Upgamma)_\mC = \Aut(\Upgamma)_{[\mC]}$.  
\end{itemize} 

Note that (AB1) is very natural: many interesting subgeometries of many point-line geometries meet this property. {\em Famous example}: let $\mC$ be a line plus 
all its points in any projective plane $\Upgamma$: then all automorphisms of the affine plane $\Upgamma \setminus \mC$ extend uniquely to automorphisms of 
the projective completion $\Upgamma$. Clearly AB-sets survive the procedure described above if the complements are confined: let $\beta$ be an automorphism of $\Upgamma \setminus \mC$; then by (AB1) it extends to an automorphism $\overline{\gamma}$ of $\Upgamma$. Since $\overline{\gamma}$ stabilizes $\Upgamma \setminus \mC$, it also stabilizes $\mC$, so by (AB3), it fixes $\mC$ elementwise. By (AB2), $\overline{\gamma}$ is trivial, so $\gamma$ also is. So what we want to construct infinite rigid projective planes, is AB-sets for which the complement is confined.

\begin{remark}{\rm 
We strongly suspect that (general, minimal and maximal) AB-sets (in various types of geometries) are interesting in their own right to study. Note that they are subgeometries which see the geometry of the ambient geometry, but \ul{don't see its symmetry}.}
\end{remark}

\subsection{Examples}

There are may interesting AB-sets (with confined complement), starting from various interesting geometries. We provide some classes of examples, starting from Kantor's planes (cf. subsection \ref{Kanpl}). 

\subsection*{In Kantor planes}

Let $\mP$ be a projective Kantor plane of order $N = 2^n$ with translation line $U$. Let $\mC = (\{ u \}, \emptyset, \I)$ be the subgeometry defined 
by the point $u$ which is not incident with $U$. Then obviously all automorphisms of $\mP \setminus \{ u \}$ extend to automorphisms of $\mP$, but $\Aut(\mP)_{\{ u \}} = \{ \id \}$ since $\Aut(\mP)$ acts sharply transitively on the points not incident with $U$. Since $\mP \setminus \mC$ is confined, 
it follows that the free completion of $\mP \setminus \mC$ is a rigid projective plane. Many other AB-sets can obviously be found in Kantor planes.\\

Note that in general (say, in geometries with more symmetry), AB-sets usually are much more complicated to find, and will need a good understanding of the fixed elements structures of elements in automorphism groups. 

 \subsection*{From rigid ovals}
 
 Let $\mP$ be a finite projective plane of order $N \geq 3$, and let $\mO$ be a rigid oval in $\mP$ (a set of $N + 1$ points no three of which are collinear, for which $\Aut(\mP)_\mO$ induces the identity on $\mO$). Then $\mO$ is an AB-set in $\mP$ and $\mP \setminus \mO$ is confined. For, let $\gamma$ be an automorphism of $\mP \setminus \mO$. We only have to show that $\alpha$ naturally defines an action on $\mO$. So let $U, V$ be two different lines which meet in a point $x \in \mO$; if $U^\gamma$ and $V^\gamma$ do not meet in a point of $\mO$, then 
 $U^\gamma \cap V^\gamma = \{ v \}$ with $v$ a point in $\mP \setminus \mO$. But then $v^{\gamma^{-1}}$ must be a point of $U \cap V$, while it is not a point of $\mO$, contradiction. 
 It follows that the action on $\mO$ is well defined, and so $\gamma$ extends to an automorphism $\overline{\gamma}$ of $\mP$. Finally, let $\alpha \in \Aut(\mP)_\mO$. As $\mO$ is rigid, it fixes every point of $\mO$. Now let $u$ be a point of $\mP \setminus \mO$, and let $U, V$ be two different lines on $u$ which meet $\mO$. If $U$ or $V$ is a tangent line 
 of $\mO$ at a point $w$ of $\mO$, then obviously it is fixed (tangent lines are unique); if $U$ or $V$ meet $\mO$ in distinct points, then obviously it is also fixed. So $u$ is also fixed, and $\alpha$ is the identity.  
 It follows that $\mO$ is an AB-set.

 The principle used here is very interesting: starting from a simple rigid structure, one constructs a more complex one. Of course, rigid ovals form just one set of rigid objects which yield AB-sets and hence rigid planes. One could also have started from rigid geometric hyperplanes in finite generalized quadrangles and use Pralle's theory of affine quadrangles \cite{Pralle} to show that these form AB-sets. 
 
 \subsection*{Rigid towers}
 
 As soon as one has constructed one rigid plane $\Upgamma_0$, it is easy to construct a countable ``tower'' of rigid planes: take away one point $x_0$, and freely construct the rigid plane $\Upgamma_1 := \overline{\Upgamma \setminus \{x\}}$; now take away one point $x_1$ in $\Upgamma_1$ to construct $\Upgamma_2 := \overline{\Upgamma_1 \setminus \{ x_1 \}}$, etc.:
 \begin{equation}
 \Upgamma_0 \ \subset \ \Upgamma_1\ \subset \ \Upgamma_2\ \ \subset\ \cdots 
 \end{equation}

\end{document}